\catcode`\@=11

\font\tenmsa=msam10 \font\sevenmsa=msam7 \font\fivemsa=msam5
\font\tenmsb=msbm10
\font\sevenmsb=msbm7 \font\fivemsb=msbm5 \newfam\msafam
\newfam\msbfam
\textfont\msafam=\tenmsa \scriptfont\msafam=\sevenmsa
\scriptscriptfont\msafam=\fivemsa \textfont\msbfam=\tenmsb
\scriptfont\msbfam=\sevenmsb \scriptscriptfont\msbfam=\fivemsb

\def\hexnumber@#1{\ifnum#1<10 \number#1\else \ifnum#1=10
A\else\ifnum#1=11
 B\else\ifnum#1=12 C\else \ifnum#1=13 D\else\ifnum#1=14
E\else\ifnum#1=15
 F\fi\fi\fi\fi\fi\fi\fi}

\def\msa@{\hexnumber@\msafam} \def\msb@{\hexnumber@\msbfam}
\mathchardef\boxdot="2\msa@00 \mathchardef\boxplus="2\msa@01
\mathchardef\boxtimes="2\msa@02 \mathchardef\square="0\msa@03
\mathchardef\blacksquare="0\msa@04 \mathchardef\centerdot="2\msa@05
\mathchardef\lozenge="0\msa@06 \mathchardef\blacklozenge="0\msa@07
\mathchardef\circlearrowright="3\msa@08
\mathchardef\circlearrowleft="3\msa@09
\mathchardef\rightleftharpoons="3\msa@0A
\mathchardef\leftrightharpoons="3\msa@0B
\mathchardef\boxminus="2\msa@0C
\mathchardef\Vdash="3\msa@0D \mathchardef\Vvdash="3\msa@0E
\mathchardef\vDash="3\msa@0F \mathchardef\twoheadrightarrow="3\msa@10
\mathchardef\twoheadleftarrow="3\msa@11
\mathchardef\leftleftarrows="3\msa@12
\mathchardef\rightrightarrows="3\msa@13
\mathchardef\upuparrows="3\msa@14
\mathchardef\downdownarrows="3\msa@15
\mathchardef\upharpoonright="3\msa@16

\mathchardef\downharpoonright="3\msa@17
\mathchardef\upharpoonleft="3\msa@18
\mathchardef\downharpoonleft="3\msa@19
\mathchardef\rightarrowtail="3\msa@1A
\mathchardef\leftarrowtail="3\msa@1B
\mathchardef\leftrightarrows="3\msa@1C
\mathchardef\rightleftarrows="3\msa@1D
\mathchardef\Lsh="3\msa@1E \mathchardef\Rsh="3\msa@1F
\mathchardef\rightsquigarrow="3\msa@20
\mathchardef\leftrightsquigarrow="3\msa@21
\mathchardef\looparrowleft="3\msa@22
\mathchardef\looparrowright="3\msa@23 \mathchardef\circeq="3\msa@24
\mathchardef\succsim="3\msa@25 \mathchardef\gtrsim="3\msa@26
\mathchardef\gtrapprox="3\msa@27 \mathchardef\multimap="3\msa@28
\mathchardef\therefore="3\msa@29 \mathchardef\because="3\msa@2A
\mathchardef\doteqdot="3\msa@2B 
\mathchardef\traceiangleq="3\msa@2C \mathchardef\precsim="3\msa@2D
\mathchardef\lesssim="3\msa@2E \mathchardef\lessapprox="3\msa@2F
\mathchardef\eqslantless="3\msa@30 \mathchardef\eqslantgtr="3\msa@31
\mathchardef\curlyeqprec="3\msa@32 \mathchardef\curlyeqsucc="3\msa@33
\mathchardef\preccurlyeq="3\msa@34 \mathchardef\leqq="3\msa@35
\mathchardef\leqslant="3\msa@36 \mathchardef\lessgtr="3\msa@37
\mathchardef\backprime="0\msa@38 \mathchardef\risingdotseq="3\msa@3A
\mathchardef\fallingdotseq="3\msa@3B
\mathchardef\succcurlyeq="3\msa@3C
\mathchardef\geqq="3\msa@3D \mathchardef\geqslant="3\msa@3E
\mathchardef\gtrless="3\msa@3F \mathchardef\sqsubset="3\msa@40
\mathchardef\sqsupset="3\msa@41
\mathchardef\trianglerighteq="3\msa@44
\mathchardef\trianglelefteq="3\msa@45 \mathchardef\bigstar="0\msa@46
\mathchardef\between="3\msa@47
\mathchardef\blacktriangledown="0\msa@48
\mathchardef\blacktriangleright="3\msa@49
\mathchardef\blacktriangleleft="3\msa@4A
\mathchardef\blacktriangle="0\msa@4E
\mathchardef\triangledown="0\msa@4F
\mathchardef\eqcirc="3\msa@50 \mathchardef\lesseqgtr="3\msa@51
\mathchardef\gtreqless="3\msa@52 \mathchardef\lesseqqgtr="3\msa@53
\mathchardef\gtreqqless="3\msa@54 \mathchardef\Rrightarrow="3\msa@56
\mathchardef\Lleftarrow="3\msa@57 \mathchardef\veebar="2\msa@59
\mathchardef\barwedge="2\msa@5A \mathchardef\doublebarwedge="2\msa@5B
\mathchardef\angle="0\msa@5C \mathchardef\measuredangle="0\msa@5D
\mathchardef\sphericalangle="0\msa@5E
\mathchardef\varpropto="3\msa@5F
\mathchardef\smallsmile="3\msa@60 \mathchardef\smallfrown="3\msa@61
\mathchardef\Subset="3\msa@62 \mathchardef\Supset="3\msa@63
\mathchardef\Cup="2\msa@64 
\mathchardef\Cap="2\msa@65
 \mathchardef\curlywedge="2\msa@66
\mathchardef\curlyvee="2\msa@67 \mathchardef\leftthreetimes="2\msa@68
\mathchardef\rightthreetimes="2\msa@69
\mathchardef\subseteqq="3\msa@6A
\mathchardef\supseteqq="3\msa@6B \mathchardef\bumpeq="3\msa@6C
\mathchardef\Bumpeq="3\msa@6D \mathchardef\lll="3\msa@6E

\mathchardef\ggg="3\msa@6F 
\mathchardef\circledS="0\msa@73
\mathchardef\pitchfork="3\msa@74 \mathchardef\dotplus="2\msa@75
\mathchardef\backsim="3\msa@76 \mathchardef\backsimeq="3\msa@77
\mathchardef\complement="0\msa@7B \mathchardef\intercal="2\msa@7C
\mathchardef\circledcirc="2\msa@7D \mathchardef\circledast="2\msa@7E
\mathchardef\circleddash="2\msa@7F
\def\ulcorner{\delimiter"4\msa@70\msa@70 }
\def\urcorner{\delimiter"5\msa@71\msa@71 }
\def\llcorner{\delimiter"4\msa@78\msa@78 }
\def\lrcorner{\delimiter"5\msa@79\msa@79 }
\def\yen{\mathhexbox\msa@55 }
\def\checkmark{\mathhexbox\msa@58 } \def\circledR{\mathhexbox\msa@72
}
\def\maltese{\mathhexbox\msa@7A } \mathchardef\lvertneqq="3\msb@00
\mathchardef\gvertneqq="3\msb@01 \mathchardef\nleq="3\msb@02
\mathchardef\ngeq="3\msb@03 \mathchardef\nless="3\msb@04
\mathchardef\ngtr="3\msb@05 \mathchardef\nprec="3\msb@06
\mathchardef\nsucc="3\msb@07 \mathchardef\lneqq="3\msb@08
\mathchardef\gneqq="3\msb@09 \mathchardef\nleqslant="3\msb@0A
\mathchardef\ngeqslant="3\msb@0B \mathchardef\lneq="3\msb@0C
\mathchardef\gneq="3\msb@0D \mathchardef\npreceq="3\msb@0E
\mathchardef\nsucceq="3\msb@0F \mathchardef\precnsim="3\msb@10
\mathchardef\succnsim="3\msb@11 \mathchardef\lnsim="3\msb@12
\mathchardef\gnsim="3\msb@13 \mathchardef\nleqq="3\msb@14
\mathchardef\ngeqq="3\msb@15 \mathchardef\precneqq="3\msb@16
\mathchardef\succneqq="3\msb@17 \mathchardef\precnapprox="3\msb@18
\mathchardef\succnapprox="3\msb@19 \mathchardef\lnapprox="3\msb@1A
\mathchardef\gnapprox="3\msb@1B \mathchardef\nsim="3\msb@1C
\mathchardef\napprox="3\msb@1D
\mathchardef\nsubseteqq="3\msb@22
\mathchardef\nsupseteqq="3\msb@23 \mathchardef\subsetneqq="3\msb@24
\mathchardef\supsetneqq="3\msb@25
\mathchardef\subsetneq="3\msb@28
\mathchardef\supsetneq="3\msb@29 \mathchardef\nsubseteq="3\msb@2A
\mathchardef\nsupseteq="3\msb@2B \mathchardef\nparallel="3\msb@2C
\mathchardef\nmid="3\msb@2D \mathchardef\nshortmid="3\msb@2E
\mathchardef\nshortparallel="3\msb@2F \mathchardef\nvdash="3\msb@30
\mathchardef\nVdash="3\msb@31 \mathchardef\nvDash="3\msb@32
\mathchardef\nVDash="3\msb@33 \mathchardef\ntrianglerighteq="3\msb@34
\mathchardef\ntrianglelefteq="3\msb@35
\mathchardef\ntriangleleft="3\msb@36
\mathchardef\ntriangleright="3\msb@37
\mathchardef\nleftarrow="3\msb@38
\mathchardef\nrightarrow="3\msb@39 \mathchardef\nLeftarrow="3\msb@3A
\mathchardef\nRightarrow="3\msb@3B
\mathchardef\nLeftrightarrow="3\msb@3C
\mathchardef\nleftrightarrow="3\msb@3D
\mathchardef\divideontimes="2\msb@3E
\mathchardef\varnothing="0\msb@3F \mathchardef\nexists="0\msb@40
\mathchardef\mho="0\msb@66 \mathchardef\thorn="0\msb@67
\mathchardef\beth="0\msb@69 \mathchardef\gimel="0\msb@6A
\mathchardef\daleth="0\msb@6B \mathchardef\lessdot="3\msb@6C
\mathchardef\gtrdot="3\msb@6D \mathchardef\ltimes="2\msb@6E
\mathchardef\rtimes="2\msb@6F \mathchardef\shortmid="3\msb@70
\mathchardef\shortparallel="3\msb@71
\mathchardef\smallsetminus="2\msb@72
\mathchardef\thicksim="3\msb@73 \mathchardef\thickapprox="3\msb@74
\mathchardef\approxeq="3\msb@75 \mathchardef\succapprox="3\msb@76
\mathchardef\precapprox="3\msb@77
\mathchardef\curvearrowleft="3\msb@78
\mathchardef\curvearrowright="3\msb@79 \mathchardef\digamma="0\msb@7A
\mathchardef\varkappa="0\msb@7B \mathchardef\hslash="0\msb@7D
\mathchardef\hbar="0\msb@7E \mathchardef\backepsilon="3\msb@7F
\def\Bbb{\ifmmode\let\next\Bbb@\else
\def\next{\errmessage{Use \string\Bbb\space only in math
mode}}\fi\next}
\def\Bbb@#1{{\Bbb@@{#1}}} \def\Bbb@@#1{\fam\msbfam#1}

\catcode`\@=\active

\def\del{\partial}

\def\CR{\hbox{{$\cal R$}}}

\def\cg{\hbox{{\sl g}}} 

\def\lform{\hbox{$\sqcup$}\llap{\hbox{$\sqcap$}}}

\def\h{{{1\over2}}}

\def\R{{\Bbb R}}
\def\C{{\Bbb C}}

\def\N{{\Bbb N}}

\def\eps{{\epsilon}}

\def\rbiprod{{\cdot\kern-.33em\triangleright\!\!\!<}}
\def\lbiprod{{>\!\!\!\triangleleft\kern-.33em\cdot\, }}

\def\tens{\mathop{\otimes}}

\def\la{{\triangleright}}

\def\extd{{{\rm d}}}

\def\isom{{\cong}}

\def\Hom{{\rm Hom}}

\def\Ad{{\rm Ad}}

\def\id{{\rm id}}

\def\<{\langle}
\def\>{\rangle}

\def\equad{\kern -1.7em}

\def\nquad{{\!\!\!\!\!\!}}

\def\eqn#1#2{\begin{equation}#2\label{#1}\end{equation}}

\def\o{{}_{\scriptscriptstyle(1)}}
\def\t{{}_{\scriptscriptstyle(2)}}
\def\th{{}_{\scriptscriptstyle(3)}}

\def\und#1{{\underline {#1}}}

\def\text#1{\mbox{\rm #1}}
\def\note#1{}

\def\blacksquare{{\lform}}
\def\frac#1#2{{{#1\over#2}}}

\def\proof{\goodbreak\noindent{\bf Proof\quad}}

\def\endproof{{\ $\lform$}\bigskip }

\def\align#1{\begin{eqnarray*}#1\end{eqnarray*}}

\def\cmath#1{\[\begin{array}{c} #1 \end{array}\]}
\def\ceqn#1#2{\begin{equation}\label{#1}
\begin{array}{c}#2\end{array}\end{equation}}

\def\vect{{\bf t}}
\def\vecx{{\bf x}}\def\vecp{{\bf p}}

\documentstyle[11pt]{article}
\newtheorem{lemma}{Lemma}[section]
\newtheorem{propos}[lemma]{Proposition}

\begin{document}
\begin{center} {\LARGE ADVANCES IN QUANTUM AND BRAIDED GEOMETRY}\\  
 \baselineskip 13pt{\ }
{\ }\\
S. Majid\footnote{Royal Society University Research Fellow and Fellow
of
Pembroke College, Cambridge}\\
{\ }\\
Department of Mathematics, Harvard University\\
Science Center, Cambridge MA 02138, USA\footnote{During the calendar
years 1995
+ 1996}\\
+\\
Department of Applied Mathematics \& Theoretical Physics\\
University of Cambridge, Cambridge CB3 9EW\\
\end{center}

\begin{center}
November 1996
\end{center}
\vspace{10pt}
\begin{quote}\baselineskip 13pt 
We demonstrate our recent general results on the Casimir construction
and moduli space of all bicovariant calculi by means of some detailed
examples, including finite-difference and 2-jet cacluli on $\R^n$ and
full details of the Casimir construction of the 4D calculus of
$SU_q(2)$. We likewise demonstrate our previous general constructions
with T. Brzezinski of quantum group gauge theory with examples of
such nonuniversal differential calculi on spacetime. We outline a
notion of quantum homotopy of a quantum space. We indicate a possible
application to classical integrable systems.
\end{quote}
\baselineskip 18pt

\section{Introduction}

There has been a lot of progress in recent years towards developing
some form
of q-geometry appropriate to q-deformed spaces. There are two main
strands; the approach based on braided categories and braid
statistics (braided
geometry) due to the author, see
\cite{Ma:book}\cite{Ma:diag}\cite{BrzMa:coa},
and the more conventional approach within  noncommutative geometry
based on
abstract differential calculus for quantum groups\cite{Wor:dif}. True
differential geometry in the latter setting, meaning connections,
gauge theory,
q-monopoles, etc. (with a quantum group fiber and quantum space base)
is due to
T. Brzezinski and the author\cite{BrzMa:gau}. By now it has
subsequently been
further developed by several authors. We emphasise the noncommutative
version
but also mention progress in the braided version and beyond it.

Until now,  emphasis has been on attempts to find q-deformations of
classical
calculi etc., or their  variants\cite{SchSch:bic}, to the extent that
this is
possible for each quantum group. The full classification problem or
construction of the `moduli space' of {\em all} possible bicovariant
differential calculi on general classes of quantum groups has
remained largely
unexplored. Some progress in the this `moduli space' direction was
made
recently in \cite{Ma:cla}, which we recall briefly in Section~2. The
main
concept is the  notion of a {\em coirreducible calculus}, without
which one
cannot begin a classification. We also recall from \cite{Ma:cla} a
new `direct'
construction of the quantum tangent spaces of bicovariant calculi
from Casimir
elements of the quantum enveloping algebra. We then provide some
detailed
examples of the classification for $\R^n$, and the Casimir
construction on
$SU_q(2)$, which are the new results beyond \cite{Ma:cla} of this
section. The
Casimir construction works uniformly in the quantum and classical
cases; we see
now why the dimension jumps when $q\ne 1$.

In Section~3 we give some detailed computations of quantum group
gauge theory
using these general nonuniversal calculi applied to the general
framework of
\cite{BrzMa:gau}. It is not hard to come up with some of these
examples by
hand, but we fit them now into the uniform general quantum group
gauge theory.
We also indicate an application of quantum group gauge theory even
with trivial
bundles, to define the quantum fundamental group of an algebra.

Although there are some modest new results, we aim here at a fresh
self-contained treatment of the geometrical theory, as a readable
introduction
to
\cite{BrzMa:gau}\cite{BrzMa:coa}\cite{Ma:diag}\cite{Ma:cla}
\cite{Haj:str}.
Some open problems will be emphasised as well. For basic results and
notations
for quantum groups themselves, see \cite{Ma:book}.

\section{Moduli of Quantum Tangent Spaces}

In this section we look at the first layer of geometry, the choice of
{\em
differentiable structure}. We are interested in our noncommutative
algebra $A$
being like the coordinate ring on a general manifold, i.e. not just a
quantum
group, though this is the case best understood at the present time.
As usual,
we do this through the notion of `tangent space' or (equivalently)
`1-forms' $\Omega^1$. By specifying one of these, we specify in
effect the
differentiable structure (which is the choice of diffeomorphism class
of a
manifold for a given topological class). The algebra by itself should
be
thought of as topological object,  needing more structure to specify
differentiation. What is interesting is that whereas  the additional
structure
is fairly unique classically (at least when we have left and right
invariance
under a group structure), one really should not expect this in the
noncommutative case, as we explain in Section~2.1. I.e. we must think
much more
in terms of the moduli space of calculi than we do classically. It is
a new
degree of freedom in quantum geometry, which we can develop field
equations
for, integrate over in q-quantum gravity, etc.

\subsection{Nonuniqueness of bicovariant calculi}

The fundamental reason for the above-mentioned non-uniqueness is as
follows. In
terms of 1-forms, we need to specify a vector space $\Omega^1$ on
which $A$
acts from the left and the right (classically, one may multiply a
1-form by a
function from the left or right). In addition, we need a linear map
$\extd:A\to
\Omega^1$ which turns `functions' into  `1-forms'. The natural axioms
are:
\begin{enumerate}
\item $\Omega^1$  an $A$-bimodule
\item The Leibniz rule $\extd(ab)=(\extd a).b + a.(\extd b)$ for all
$a,b\in A$
\item Surjectivity in the form $\Omega^1={\rm span}\{a\extd b\}$
\end{enumerate}
This is called the choice of {\em first order differential calculus}.
It is
more or less the minimal structure needed for some kind of
differential
geometry. Note that when $A$ is not commutative, it makes no sense to
assume
that
\eqn{comcalc}{a.(\extd b)=(\extd b).a}
as we would assume classically. This is because such an assumption
and the
Leibniz rule would imply $\extd(ab)=\extd(ba)$, which is not
reasonable when
our algebras are non-commutative. This is the reason that we require
$\Omega^1$
to be an $A$-bimodule. On the other hand, the assumption
(\ref{comcalc}) of a
`commutative calculus' is highly restrictive classically, and when we
relax it
we have a much larger range of possibilities than in the classical
case.

When $A$ is a Hopf algebra it is natural to focus on differential
calculi which
are covariant under left or right translations on the `group'.
Classically,
this invariance condition fixes the differentiable structure uniquely
as that
obtained by translation from the differentiable structure at the
identity. One
usually identifies the tangent space at the identity with the Lie
algebra $\cg$
and extends it to $G$ via left-invariant vector fields $\tilde{\xi}$
for
$\xi\in G$. These then transform under right translation via the
adjoint
action. Likewise, at the level of `1-forms'. One may impose similar
conditions
in the quantum case, that $\Omega^1$ is also an $A^*$-bimodule (or
more
precisely an $A$-bicomodule) in a manner compatible with the
$A$-bimodule and
$\extd$ structures. This is called a {\em bicovariant differential
calculus}.
Being the most restrictive, we might hope to have uniqueness at least
in this
case.

We recall first that every quantum group has a counit $\eps:A\to\C$
which,
classically, denotes evaluation at the group identity. The space
$\ker\eps\subset A$ has a natural action of $A$ by multiplication
form the left
and of $A^{*\rm op}$ by the quantum coadjoint action $\Ad^*$. The two
actions
form an action of the quantum double $D(A)$. Then, after a little
analysis, one
finds cf\cite{Wor:dif}:

\begin{propos} A bicovariant $\Omega^1$ must be of the form
$\Omega^1=V\tens A$
\cmath{ a.(v\tens b)=a\o\la v\tens a\t b,\quad
(v\tens b).a=v\tens ba,\quad \extd a=(\pi\tens \id)(1\tens a-a\o\tens
a\t)}
for some vector space $V$ which is a quotient of $\ker\eps$ as a
$D(A)$-module,
i.e. $V=\ker\eps/M$ where $M\subset \ker\eps$ is some
$D(A)$-submodule. Here
$\pi:\ker\eps\to V$ denotes the canonical projection and $\la$
denotes
left multiplication by $A$ projected down to $V$.
\end{propos}

The vector space $V$ is called the space of left-invariant 1-forms,
and its
dimension is called the dimension of the differential calculus. For
$A=SU_q(2)$
the lowest dimension  bicovariant calculus when $q\ne 1$ has
dimension 4 (not
3!) and was found in \cite{Wor:dif}. Woronowicz found in fact two
calculi of
dimension 4 given by a similar construction. More generally, for the
ABCD
series of semisimple Lie groups there is a natural R-matrix
construction\cite{Jur:dif} for a bicovariant calculus of dimension
$n^2$ (where
$G\subset M_n$), and some variants of it which are classified
case-by-case in
\cite{SchSch:bic}.

It may seem therefore that (up to some variants related to a choice
of roots)
there is one natural bicovariant calculus for each of the standard
quantum
groups $G_q$. This is misleading, however. It was shown in
\cite{BrzMa:bic}
that there is at least a whole `algebra' of calculi associated to the
non-trivial elements $\alpha$ of the algebra $Z^*(A)$ of
$\Ad^*$-invariant
elements. Classically, these are the class functions on the group $G$
(the
functions which are constant on each conjugacy class). Explicitly,
\eqn{alphaid}{\alpha\in Z^*(A)=\{a\in A|  a\o(Sa\th)\tens a\t=1\tens
a\},\quad
V =\ker\eps/\ker\eps(\alpha-(\eps(\alpha)+1))}
defines the calculus, where $(\id\tens\Delta)\circ\Delta a=a\o\tens
a\t\tens
a\th$ is the notation. Moreover,  $\extd$ then has the inner
form\cite{BrzMa:bic}
\eqn{inner}{ \extd a=a.\omega_\alpha-\omega_\alpha.a,\quad
\omega_\alpha=(\extd\alpha\o)S\alpha\t\in\Omega^1,}
where $S$ is the antipode or `linearised inversion' operator on $A$.
This shows
that the entire construction is not possible classically, since
(\ref{comcalc})
would force $\extd a=0$. For the standard quantum groups $G_q$, it is
also
shown in \cite{BrzMa:bic} that $Z^*(G_q)$ is commutative with
rank-$\cg$
algebraically independent generators $\{\alpha_i\}$ (this is proved
by studying
the braided group version of $G_q$). Only one of these, the
$q$-trace, is
relevant to the natural $n^2$-dimensional calculus of \cite{Jur:dif},
which
appears after making further quotients of $V$\cite{BrzMa:bic}.

For example, for the quantum group $SU_q(2)$ with the standard matrix
of
generators $\vect=\pmatrix{a&b\cr c&d}$, we take
\eqn{alphasu2}{ \alpha={q a+ q^{-1}d\over(q-q^{-2})(q-1)}.}
Here $\extd$ in (\ref{inner}) has a classical remnant as $q\to 1$
because the
normalisation of $\alpha$ goes to $\infty$, yielding a certain
bicovariant differential calculus on $SU_2$.

\subsection{Construction of invariant quantum tangent spaces}

Even the family of calculi associated to the algebra $Z^*(A)$ does
not exhaust
the bicovariant calculi on a typical quantum group. To get some feel
for the
entire moduli space, i.e. to classify all calculi, it is convenient
to
concentrate on the {\em invariant quantum tangent space}
associated to every bicovariant calculus. This is the `tangent space
at the
identity' $L=V^*$ and can be characterised as follows. First of all,
note that
$A^*$ (defined in a suitable way) will generally also be a Hopf
algebra, with
its own counit $\eps$. This is the self-duality of the axioms of a
Hopf
algebra. Moreover, $\ker\eps\subset A^*$ is acted upon by the quantum
double
$D(A^*)$ by
\eqn{Laction}{ x\la y=\Ad_x(y)=x\o y Sx\t,\quad a\la y=\<a,y\o\>y\t -
1\<a,y\>,\quad \forall y\in\ker\eps\subset A^*}
for the action of $a\in A$ and $x\in A^*$. Here $\Ad$ is the quantum
adjoint
action, while the action of $A$ is rather strange. Classically,
$A^*=U(\cg)$
and $\ker\eps$ denotes any product of Lie algebra elements. The
quantum adjoint
action becomes the usual coadjoint action restricted to $\ker\eps$,
and the
action of $A=\C(G)$ becomes to lowest degree
\[ a\la \xi=a(e)\xi,\quad
a\la(\xi\eta)=a(e)\xi\eta+\<a,\xi\>\eta+\<a,\eta\>\xi\]
etc for $\xi,\eta\in \cg$. The action on higher degree elements of
$U(\cg)$ is
more complicated, as determined by the coproduct $\Delta$.

\begin{propos} The possible invariant quantum tangent spaces are
precisely
the subspaces $L\subset\ker\eps\subset A^*$ which are
subrepresentations of the
quantum double $D(A^*)$ acting on $\ker\eps$.
Explicitly, this means $L$ such that:
\eqn{Lie}{ L\subset\ker\eps\subset A^*,\quad \Ad_x(L)\subset L,\quad
(\Delta
-\id\tens 1)(L)\subset  A^*\tens L}
for all $x\in A^*$.
\end{propos}
The explicit version (\ref{Lie}) works (more or less equivalently)
with
everything in terms of the Hopf algebra structure of $A^*$ rather
than
mentioning the quantum double.  Also note that although these quantum
tangent
spaces $L$ are sometimes called `quantum Lie algebras'\cite{Wor:dif},
they do
not usually have enough structure to qualify for a Lie algebra of any
kind. In
nice cases
where they do have enough structure, they are more properly
formulated as
{\em braided-Lie
algebras}\cite{Ma:lie}. Moreover, for any Hopf algebra $A$, the
associated
differentiation operators $\del_x:A\to A$ defined by
$\del_x(a)=(\<x,\
\>\tens\id)\circ\extd$ obey a {\em braided Leibniz rule}\cite{Ma:cla}
\eqn{bleib}{\del_x(ab)=(\del_x a)b +a_i\del_{x^i} b}
for all $a,b\in A$ and $x\in L$. Here $\Psi^{-1}(a\tens x)\equiv
x^i\tens a_i$
is our notation (summation understood) and $\Psi$ is the braiding of
$L$ with
$A$ as spaces on which the quantum double $D(A^*)$ acts. The elements
of $L$
define in this way the `braided left-invariant vector fields'
associated to a
bicovariant calculus.

This quantum tangent space point of view is developed in detail in
\cite{Ma:cla}, where some general classification theorems are then
obtained.
Firstly, one should introduce the notion of {\em coirreducible
calculus}, which
we define as corresponding to an irreducible quantum tangent space.
It turns
out that the classical case $A=\C(G)$ is rather singular and is
actually the
hardest; we say more about it in Section~2.3. The finite group
classical case
is easier and it is known already\cite{BMD:non} that natural
bicovariant
calculi are obtained from nontrivial conjugacy classes in $G$. It is
proven in
\cite{Ma:cla}
that these are in fact coirreducible, and that coirreducible calculi
on
$A=\C(G)$ are in
1-1 correspondence with nontrivial conjugacy classes. The span of the
conjugacy
class is
$L$. Also, taking the view of noncommutative geometry that
the group algebra $A=\C G$ for a finite group can be viewed as `like'
$\C(\hat
G)$ (here $\hat G$ need not exist, however), one finds\cite{Ma:cla}
that the
coirreducible calculi in this case are classified by pairs consisting
of an
irreducible representations $\rho$ and a continuous parameter in the
projective space $\C P^{\dim\rho-1}$. They have dimension $(\dim
\rho)^2$.
Finally,
for a strict semisimple quantum group (with a universal R-matrix
$\CR$ obeying
some further strict conditions) one can again classify the
coirreducible
calculi\cite{Ma:cla}; they turn out again to correspond to the
irreducible
representations $\rho$ of $\cg$ and have dimension $(\dim \rho)^2$.
This
applies essentially to the standard deformations $G_q$ for generic
$q$, up to
some variants. In the case of $SU_q(2)$, one knows that the quantum
double of
$U_q(su_2)$ can be identified when $q\ne 1$ with some version of the
$q$-Lorentz group. Hence its possible invariant quantum tangent
spaces are the
$q$-deformation (and projection to $\ker\eps$) of the
subrepresentations of the
particular representation consisting of the Lorentz group acting on
the space
of functions on $SU_2$ by left and right
group translation. The lowest possible coirreducible calculus is
therefore
4-dimensional (the action on Minkowski space). This is the tensor
square of
the spin $1/2$ representation of $SU_2$. The full
analysis\cite{Ma:cla}
shows that generically (up to some uniform choices), the tensor
square of each
spin $j>0$ irreducible representation of $SU_2$ occurs, just once. So
there is
a
natural 9-dimensional spin 1 calculus, a 16 dimensional spin 3/2
calculus, etc.

Hence the classical $q=1$ theory is more like the finite group case,
while the
$q\ne 1$ theory is more like the finite group-dual case. In all
cases, there is
an entire moduli space of calculi. One may endow the moduli space
with a
topology, and introduce natural operations on it. One of them, in
\cite{Ma:cla}, associates to every bicovariant calculus a
`dual' or `mirror' calculus of a quite different form but on the same
quantum
group.

Also introduced in \cite{Ma:cla} is a natural construction for
bicovariant
quantum tangent spaces which is `dual' to the $\alpha$ construction
in
Section~2.1. It associates a quantum tangent space to any non-trivial
element
in the centre $Z(A^*)$. This is the algebra of elements of $A^*$
fixed under
the quantum adjoint action.

\begin{propos}\cite{Ma:cla} For any $c\in Z(A^*)$,
\eqn{cascalc}{ L={\rm span}\{x_a=\<a,c\o\>c\t -\<a,c\>|\ a\in\
\ker\eps \subset
A\}}
defines an invariant quantum tangent space.
\end{propos}

We now compute the case $A=SU_q(2)$ and $A^*=U_q(su_2)$ in detail. We
take the
latter in its standard form with generators $q^{H\over 2},X_\pm$,
relations,
coproduct and counit
\ceqn{Uqsu2}{ \nquad q^{H\over 2}X_\pm q^{-{H\over 2}}=q^{\pm
1}X_\pm,\
[X_+,X_-]={q^H-q^{-H}\over q-q^{-1}},\  \Delta q^{H\over 2}=q^{H\over
2}\tens
q^{H\over 2},\quad\eps\, q^{H\over 2}=1\\
 \Delta X_\pm=X_\pm\tens q^{H\over 2}+q^{-{H\over 2}}\tens
X_\pm,\quad \eps\,
X_\pm=0.}
There is a well-known quadratic $q$-Casimir which we take in a
certain
normalisation and offset by a multiple of 1:
\eqn{cassu2}{ c_q={C-(q+q^{-1})\over (q-q^{-2})(q-1)};\quad
C=q^{H-1}+q^{-H+1}+(q-q^{-1})^2 X_+X_-.}

\begin{propos} The invariant quantum tangent space of $SU_q(2)$
generated by
the $q$-Casimir $c_q$ in (\ref{cassu2}) via Proposition~2.3 is
4-dimensional
and coincides with the braided Lie algebra $\und{gl}_{2,q}$ in
\cite{Ma:lie}
when
$q\ne 1$. Explicitly, $ L={\rm span}\{x_{a-1},x_b,x_c,x_{d-1}\}$,
where
\cmath{x_{a-1}={q+1\over q-q^{-2}}(q^H-1) + (q^{-1}-1)c_q,\qquad
x_{b}={q^\h(q+1)(1-q^{-2})\over q-q^{-2}}q^{H\over 2}X_-\\
x_c= {q^\h(q+1)(1-q^{-2})\over q-q^{-2}}X_+ q^{H\over 2},\qquad
x_{d-1}={q^{-1}+1\over q^{-1}-q^2}(q^H-1)+(q-1)c_q.}
When $q=1$, the space $L$ is three-dimensional and coincides with the
Lie
algebra $su_2$.
\end{propos}
\proof The duality pairing of $U_q(su_2)$ with $SU_q(2)$ is the
fundamental
representation
\cmath{\<\vect,q^{H\over 2}\>=\pmatrix{q^\h&0\cr 0&
q^{-\h}},\quad \<\vect,X_+\>=\pmatrix{0&1\cr 0& 0},\quad
\<\vect,X_-\>=\pmatrix{0&0\cr 1& 0}.}
Using this and the identity
\align{ \Delta C\equad &&=C\tens q^H+q^{-H}\tens
C-(q+q^{-1})q^{-H}\tens q^H \\
&&\qquad\qquad +(q-q^{-1})^2(X_+q^{-{H\over
2}}\tens q^{H\over 2}X_-+q^{-{H\over 2}}X_-\tens X_+q^{H\over
2})}
which follows from (\ref{Uqsu2}), we compute the elements of $L$
corresponding
to $a-1,b,c,d-1\in\ker\eps\subset SU_q(2)$. They are clearly linearly
independent when $q\ne1$. The products of these generators with
general
elements of $SU_q(2)$ taken in the role of $a\in\ker\eps$ in
Proposition~2.3 do
not give
new elements of $L$ because this 4-dimensional $L$ is already known
to be
stable under the action of all $a\in SU_q(2)$ in (\ref{Laction}).
Hence these
elements are a
basis for $L$. This vector space coincides with the braided Lie
algebra $\und
{gl}_{2,q}$ because this is identified\cite{Ma:skl} with a
natural subspace of $U_q(su_2)$, according to
\ceqn{lieenv}{ h={q^{-1}\over q^2-1}(C-(q+q^{-1})q^H),\quad
x=q^{-{3\over
2}}q^{H\over 2}X_-\\ y=q^{-{3\over 2}}X_+q^{H\over 2},\quad
\gamma={q^{-1}\over q^2-1}(C-(q+q^{-1})).}
We recall that the braided Lie bracket in \cite{Ma:skl}\cite{Ma:lie}
is given
by the quantum adjoint action restricted to this subspace.
Explicitly,
\ceqn{bralie}{{}[h,x]=(q^{-2}+1)x=-q^{2}[x,h],\quad
[h,y]=-(q^{-2}+1)q^{-2}y=-q^{-2}[y,h]\\
{}[x,y]= q^{-2}h=-[y,x],\quad [h,h]=(1-q^{-4})h,
\quad [\gamma,\cases{h\cr x\cr y}]=(1-q^{-4})\cases{h\cr x\cr y}}
The isomorphism (\ref{lieenv}) is valid for generic $q\ne
1$. By contrast, as $q\to 1$ the elements $x_{a-1}$ and $x_{d-1}$
coincide and
$L$ becomes 3-dimensional in this special limit.
\endproof

The braided-Lie algebra $\und{gl}_{2,q}$ is known to be the  quantum
tangent
space for the usual 4-dimensional bicovariant differential calculus
on
$SU_q(2)$. The Casimir construction  in Proposition~2.3 from
\cite{Ma:cla}
provides a new
and more direct,
route. Moreover, in the moduli of all quantum tangent spaces, only
very few
calculi (probably just one in the standard cases) will be close
enough to the
classical calculus to have the additional properties needed properly
to be a
`quantum' or braided Lie
algebra. Using the q-deformed quadratic Casimir of $U_q(\cg)$ exactly
picks out
such a natural
one, as the above example demonstrates in detail. Moreover, we
achieve this
without R-matrices; it works for general simple $\cg$,
including the exceptional series.

\subsection{Quantum tangent spaces on classical groups}

In this section we specialise to the case $A=\C(G)$, where $G$ is a
classical
Lie group. More precisely, we take $A^*=U(\cg)$ where $\cg$ is the
Lie algebra
of $G$. The classification of bicovariant calculi or their quantum
tangent
spaces in Proposition~2.2 is not known in any generality:
\begin{itemize}
\item It is an {\em unsolved problem in classical Lie theory to find
all
irreducible subspaces $L\subset U(\cg)$ stable under $\Ad$ and
$\Delta-\id\tens
1$}.
\end{itemize}

Of course, many possible $L$ (not necessarily irreducible) are known.
For
example, we may take
\eqn{jetg}{L=\cg,\quad L=\cg+\cg\cg,\quad L=\cg+\cg\cg+\cg\cg\cg}
etc. as subspaces of $U(\cg)$. Only the first of these obeys
(\ref{comcalc})
but the rest are just as valid as soon as we relax this assumption.
The
construction in Proposition~2.3 works perfectly well in this
classical case
i.e. the q-deformed case in Section~2.2 has a smooth limit as $q\to
1$. When
$\cg$ is semisimple, there are rank-$\cg$ algebraically independent
Casimirs,
each with a natural bicovariant calculus.

We will focus in fact on the very simplest case, which nevertheless
demonstrates the key
points: we take $G=\R$ the additive real line. Let us note that the
possibility
and applications of non-standard differential calculi on classical
spaces such
$\R^n$ have already been emphasised in the papers of M\"uller-Hoissen
and
collaborators\cite{DimMul:non}\cite{DMS:non}, although perhaps not as
part of a
systematic theory of bicovariant calculi and quantum tangent spaces
as we
present now.

Thus, we take $A=\C[x]$ and $U(\R)=\C[p]$, i.e. functions in
1-variable with
their linear coproducts. According to Proposition~2.2, we need to
classify all
subspaces $L\subset\C[p]$ which obey $\forall f(p)\in L$,
\eqn{LR}{f(0)=0,\quad  (\Delta-\id\tens 1)f=f'\tens p + f''\tens
{p^2\over
2!}+f'''\tens {p^3\over 3!}+\cdots\in \C[p]\tens L.}
An obvious family of choices is $L={\rm span}\{p,p^2,\cdots,p^n\}$
for any
natural number $n$. For then all $f\in L$ have the property that
$f^{(m)}=0$ if
$m>n$, ensuring the required condition. One can also write the second
condition
in (\ref{LR}) as
$f(p+\lambda)-f(\lambda)\in L$ for all $\lambda$ (by evaluating
against
$\lambda$). So

\begin{itemize}
\item An invariant quantum tangent space for $\R$ means a subspace of
functions vanishing at zero and closed under all projected
translations (i.e.
translations followed by subtraction of the value at zero).
\end{itemize}

We can also approach the construction of natural $L$ through
Proposition~2.3.
All functions $c(p)$ are central, so every function determines a
calculus! This
is therefore a new degree of freedom determined by an arbitrary
function on
$\R$. For a basis of $\ker\eps\subset A$, we take $\{x,x^2,\cdots\}$.
Then,

\begin{propos} Any function $c(p)$ defines an invariant quantum
tangent space
on $\R$,
\[ L={\rm span}\{ p_n=c^{(n)}(p)-c^{(n)}(0),\ n\in\N\}.\]
The associated braided derivations $\del_{p_n}$ and their braiding
with
functions $f(x)$ are
\[ \del_{p_n}f=(c^{(n)}({\extd\over\extd x})-c^{(n)}(0))f,\quad
\Psi^{-1}(f\tens p_n)=\sum_{m=0}^\infty p_{n+m}\tens {f^{(m)}\over
m!}\]
\end{propos}

The associated space of invariant 1-forms $V$ is defined as $L^*$ and
is
computed as follows. We start with $\ker\eps\subset \C[x]$, i.e. with
$f(x)$
such that $f(0)=0$, and quotient out by the subspace of $f$ such that
\eqn{VR}{ (\del_{p_n}f)(0)=0,\quad \forall n\in\N.}
This follows from the duality $\<x^n,p^m\>=\delta_{n,m}m!$ between
$\C[x]$ and $\C[p]$. We then choose a basis $\{e_i\}$ of $L$ with
dual basis $\{f^i\}$ and define
\eqn{extdR}{\extd:A\to\Omega^1=V\tens A,\quad \extd f=f^i\tens
\del_{e_i} f.}

\begin{propos} The function $c(p)=p^{n+1}/(n+1)!$ in Proposition~2.5
yields the
$n$-dimensional calculus with $L={\rm span}\{p,p^2,\cdots , p^n\}$.
The
derivations, their braiding and exterior derivative are
$\del_{p^m}f=f^{(m)}$
and
\cmath{ \Psi^{-1}(f\tens
p^m)=\sum_{r=0}^{m-1}p^{m-r}\tens f^{(r)}{m!\over r!(m-r)!},\quad
\extd
f=\sum_{m=1}^n {\pi(x^m)\over m!}\tens f^{(m)}(x)}
where $\pi$ is the projection defined by $x^{n+1}=0$.
\end{propos}
\proof The $p_n$ generators of $L$ from Proposition~2.5 are
$p_{n-m+1}=p^m/m!$
for $1\le m\le n$ and $p_m=0$ for $m>n$. Hence we can take
$\{p,\cdots,p^n\}$
as basis of $L$. Their action from the duality pairing between
$\C[x],\C[p]$ is
by $p$ as usual
differentiation, while the braiding $\Psi^{-1}$ follows immediately
from the
formula in Proposition~2.5. The space $V$ consists of functions
vanishing at
zero modulo powers of $x$ higher than $x^n$. Hence
$\{\pi(x),\pi(x^2)/2!,\dots,\pi(x^n)/n!\}$ is a dual basis for $V$.
We then
obtain $\extd f$ from (\ref{extdR}). \endproof

The case $n=1$ recovers the usual differential calculus on $\R$,
while $n>1$
gives higher order {\em jet calculi} in which higher derivatives are
regarded
as first-order vector fields, but obeying a braided Leibniz rule
(\ref{bleib})
with nontrivial $\Psi^{-1}$. The commutativity (\ref{comcalc}) also
fails for
$n>1$. For example,  $c(p)=p^3/3!$ gives the 2-jet calculus with
1-forms
obeying
\eqn{2-jet}{ \extd f=(\extd x)f'+\omega f'';\quad
\omega\equiv\h(x\extd
x-(\extd x)x),\quad {}f\extd x -(\extd x)f=2\omega f',\quad
f\omega=\omega f,}
where we identify $\pi(x)\tens 1=\extd x$ and $\h\pi(x^2)\tens
1=\omega$ as
invariant forms.
For another example, we formally make completions of our algebras to
allow
certain powerseries. Then

\begin{propos} For all $\lambda\ne0$, the function
$c(p)=\lambda^{-2}e^{\lambda
p}$ in Proposition~2.5 yields a 1-dimensional bicovariant
differential calculus
on $\R$ with $L=\{p_1=\lambda^{-1}(e^{\lambda p}-1)\}$. The
differentiation and
braiding are
\[ \del_{p_1}f={f(x+\lambda)-f(x)\over\lambda},\quad \Psi^{-1}(f\tens
p_1)=p_1\tens f(x+\lambda).\]
\end{propos}
\proof From Proposition~2.5, we find $p_{n+1}=\lambda^n p_1$ so $L$
is
1-dimensional. The braiding then follows at once from $\Psi^{-1}$ in
Proposition~2.5 and Taylor's theorem. \endproof

Thus, the finite-difference operations introduced by Newton's tutor
Barrow {\em are exact} differentiation with but respect to a
non-standard
differential
calculus on $\R$. In fact, all coirreducible calculi on $\R$ are of
this form,
parametrised
by $\lambda$. And we only needed
to relax (\ref{comcalc}) in order to allow them, without taking the
limit
$\lambda\to0$. To see that (\ref{comcalc}) fails, we compute the
1-forms as
follows. The vector space $V$ consists of functions vanishing at $0$
quotiented
out by all $f$ such that $f(\lambda)=0$. This is essentially the
functions
divisible by $x$ modulo those divisible by $x(x-\lambda)$, which is a
1-dimensional space. As a basis of $V$ we can take  $\{\pi(x)\}$
since
$\<p_1,x\>=1$. Then
\eqn{chordform}{ \extd f= \pi(x)\tens
{f(x+\lambda)-f(x)\over\lambda}=(\extd
x){f(x+\lambda)-f(x)\over\lambda},}
where we identify $\extd x=\pi(x)\tens 1$ to obtain the second
expression.
Similarly, we find that $x\extd x-(\extd x)x=\pi(x^2)\tens
1=\lambda\pi(x)\tens
1=\lambda \extd x$, so we see that (\ref{comcalc}) fails for this
calculus.
More generally, we deduce from this that
\eqn{chordnoncom}{  f\cdot \extd x - (\extd x)\cdot f=\lambda \extd
f}
so that the calculus is of the inner form (\ref{inner}) discussed in
\cite{BrzMa:bic}.

The story is similar for $\R^n$. For example, on $\R^2$ the function
$c(\vecp)=\lambda^{-2}e^{\lambda p}+\mu^{-2}e^{\mu q}$ (where
$\vecp=(p,q)$ and
$\lambda,\mu\ne 0$) gives the 2-dimensional calculus $L={\rm
span}\{p_{1,0},p_{0,1}\}$ with
\ceqn{2-dif}{
\del_{p_{1,0}}f={f(x+\lambda,y)-f(x,y)\over\lambda},\quad
\del_{p_{0,1}}f={f(x,y+\mu)-f(x,y)\over \mu},\quad \extd f=(\extd x)
\del_{p_{1,0}} f+ (\extd y)\del_{p_{0,1}} f\\
x\extd y=(\extd y)x,\quad y\extd x=(\extd x)y,\quad x\extd x-(\extd
x)x=\lambda\extd x,\quad y\extd y-(\extd y)y=\mu\extd y.}
Of course, there are many other interesting functions $c(p)$ one
could take.
Another interesting one is $c(p)=\lambda^{-1}e^{\lambda p^2}$
yielding a
`Hermite differential calculus'. We truly have a new field in physics
determining the differential calculus on spacetime.

We have considered here the bicovariant calculi on quantum groups. It
is
straightforward
to write down the corresponding axioms for a braided group $B$ and
make a
corresponding
analysis. The role of the quantum double in the general theory is
played by
the author's double-bosonisation Hopf algebra $B\lbiprod H\rbiprod
B^{*\rm op}$
in \cite{Ma:dbos} (where $H$ is the background quantum group in the
braided
category generated by which $B$ lives).
One should think of it as the bosonisation of the braided double
$D(B)$,
although this does not
really exist by itself when the category is truly braided. Detailed
diagrammatic axioms and proofs will be presented elsewhere. But the
final
result (as presented in Prague in June 1996) is, in suitable
conventions,

\begin{propos} The possible invariant braided tangent spaces are
precisely
the subspaces $L\subset \ker\eps\subset B^*$ which are
subrepresentations under
the
double-bosonisation  $B\lbiprod H\rbiprod B^{*\rm op}$. Here $B^*$
acts by the
braided adjoint action, $B$ by evaluation against projected
translation and $H$
acts as the background quantum group, cf. the canonical action in
\cite{Ma:dbos}.
\end{propos}

This is the braided version of Proposition~2.2. On the other hand,
when we look
at the q-deformed (braided) bicovariant
differential calculi on $B=B^*=\R_q^n$, we have an entirely different
story
when
$q\ne 1$. In this case the double-bosonisation is the $q$-conformal
group\cite{Ma:conf} (and the background quantum group is the
dilaton-extended
$q$-rotation group). The possible invariant braided quantum tangent
spaces $L$
are q-deformations of the projections (to vanish at zero) of
irreducible
subrepresentations of
the conformal group acting on functions on $\R^n$. This means that
$L$ is
more or less unique in this case and is (the projection of) the
q-deformed set
of solutions of the massless Klein-Gordon equation $\lform f=0$. This
uniqueness gives us a clue to what will be involved in `gluing'
copies of
$\R_q^n$
together to define manifolds if we want to keep everything
bicovariant under
translation: infinite-dimensional tangent spaces defined by solutions
of the
massless wave equation. If we relax translational bicovariance then
other
q-deformations of the 1-forms on $\R^n$ are also possible, including
one which
extends to the whole exterior algebra and which typically has the
classical
dimensions in each degree; see\cite{Ma:varen}\cite{Ma:book}. Some
recent
categorical results about exterior algebras on braided groups are in
\cite{BesDra:dif}.

\section{Gauge Theory}

 Once we have fixed choices of differential calculi on our various
spaces, we
can do a systematic quantum group gauge theory. This formalism has
been
introduced in \cite{BrzMa:gau}, where the $q$-monopole over the
$q$-deformed
$S^2$ was constructed. We work with all `spaces' in terms of
algebras. So the
`total space' of a principal bundle is an algebra $P$. A Hopf algebra
$A$
coacts (or $A^*$ acts) on $P$ and the fixed point subalgebra
$M\subseteq P$
plays the role of the `base space' of the bundle. Among the 1-forms
on $P$ we
have the horizontal 1-forms $P\Omega^1(M)P$ and a {\em connection} is
a choice
of complement to this. One also has a theory of gauge transforms,
gauge-covariant curvature, associated vector bundles and connections
on them;
see \cite{BrzMa:gau}. The theory has been generalised to braided
groups and,
beyond, to systems where $A$ is merely a
coalgebra\cite{BrzMa:coa}\cite{Ma:diag}.

Here we limit ourselves to gauge theory for the case
of trivial bundles of the form $P=M\tens A$. In this case, all
formulae look
like more usual gauge theory. Thus, we consider
gauge fields, curvature, gauge transformations and matter fields as
\eqn{trivcon}{ \alpha\in \Omega^1(M)\tens A^*,\quad
F\in\Omega^2(M)\tens
A^*,\quad \gamma\in M\tens A^*,\quad \psi\in \Omega^n(M)\tens V}
respectively. Here $A^*$ is assumed to act in $V$. At this level, we
do not
really need a differential calculus on $A$ itself. We do not even
need $A$ to
be a Hopf algebra. The minimum is that $A$ should be a coalgebra, or
$A^*$ an
algebra. When $A$ is a Hopf algebra (or even a braided group) we
should chose a
differential calculus $\Omega^1(A)$ and define
$\Omega^1(P)=P\Omega^1(M)P\oplus
P\Omega^1(A)$. This splitting corresponds to a canonical trivial
connection on
$P=M\tens A$, while other connections are obtained from this by
adding
$\alpha$\cite{BrzMa:gau}. One can also try to restrict
$\alpha\in\Omega^1(M)\tens L$, although this requires some further
work. None of this
is needed for bare-bones gauge theory, where we work directly on the
base and
forget about the full geometrical structure of the principal bundle.

\subsection{Quantum homotopy group $\pi_1(M)$}

Let $M$ and $A^*$ be algebras. We outline the generalised gauge
theory at this
level, and discuss a first application. As explained, we need only
$\Omega^1,\Omega^2$ on the base algebra $M$. By the former, we mean a
choice of
first order differential calculus over $M$ as at the start of
Section~2.1. For
the 2-forms we need a vector space $\Omega^2$ such that:

\begin{enumerate}
\item $\Omega^2$ is an $M$-bimodule.
\item There is a surjection
$\wedge:\Omega^1\tens_M\Omega^1\to\Omega^2$ as
$M$-bimodules.
\item There is a map $\extd:\Omega^1\to \Omega^2 $ obeying
$\extd\circ\extd=0$
when composed with $\extd:M\to\Omega^1$.
\item The Leibniz rule $\extd(a.\omega.b)=(\extd a)\wedge
\omega.b+a.(\extd\omega).b-a.\omega\wedge\extd b$ for all
$\omega\in\Omega^1$
and $a,b\in M$.
\end{enumerate}

The first two items mean that $\Omega^2$ is a quotient of
$\Omega^1\tens_M\Omega^1$. The third and fourth conditions fully
specify
$\extd$ once we have chosen this quotient, according to
\eqn{comm2form}{ \extd(a\extd b)=\extd a\wedge \extd b,\quad
\extd((\extd
b)a)=-\extd b\wedge\extd a.}
In particular, we see that if (\ref{comcalc}) holds then $\wedge$ is
antisymmetric (as it is classically), but in general we cannot assume
this.
Also, one
knows that $\Omega^1$ can be built by quotienting a certain universal
calculus
on $M$ by a bimodule. This determines likewise a natural choice of
the quotient
of $\Omega^1\tens_M\Omega^1$ which defines $\wedge$. In practice,
these minimal
relations of $\Omega^2$ are simply obtained by applying $\extd$
(extended by
the Leibniz rule) to the relations of $\Omega^1$. We will use
$\Omega^2$
defined canonically from $\Omega^1$ in this way; one could consider
further
quotients too.

With these basic structures fixed, we can define the essential
features of
gauge theory. We consider gauge fields, gauge transformations etc. as
in
(\ref{trivcon}) and define gauge transformations, curvature and the
covariant
derivative by
\eqn{gaugetheory}{
\alpha^\gamma=\gamma^{-1}\alpha\gamma+\gamma^{-1}\extd\gamma,\quad
\psi^\gamma=\gamma^{-1}\psi,\quad F(\alpha)=\extd \alpha +
\alpha\wedge\alpha,\quad \nabla\psi=\extd\psi+\alpha\wedge\psi.}
Products in $A^*$ and its action on $V$ are understood, along with
the above
operations involving $M$. For example, the allowed gauge
transformations are
the invertible elements of the algebra $M\tens A^*$. One can then
verify the
fundamental lemmas of gauge theory, namely

\begin{enumerate}
\item $F(\alpha^\gamma)=\gamma^{-1}F(\alpha)\gamma$.
\item $\nabla^\gamma\psi^\gamma=(\nabla\psi)^\gamma$
\item $\nabla^2\psi=F\psi$
\end{enumerate}

One also has Bianchi identities, etc. at this level\cite{BrzMa:gau}.
For
example, to see exactly what is involved in the gauge covariance of
$F$, we
note first that
$\extd\gamma^{-1}=-\gamma^{-1}(\extd\gamma)\gamma^{-1}$ from
the Leibniz rule. Then
\align{F(\alpha^\gamma)\equad
&&=\extd(\gamma^{-1}\alpha\gamma+\gamma^{-1}\extd\gamma)
+(\gamma^{-1}\alpha\gamma+\gamma^{-1}\extd\gamma)\wedge
(\gamma^{-1}\alpha\gamma+\gamma^{-1}\extd\gamma)\\
&&=-\gamma^{-1}(\extd
\gamma)\gamma^{-1}\wedge\alpha\gamma+\gamma^{-1}(\extd\alpha)\gamma-\
gamma^{-1}
\alpha\wedge\extd\gamma-\gamma^{-1}\extd\gamma\wedge
\gamma^{-1}\extd\gamma+\gamma^{-1}\extd^2\gamma\\
&&\quad+\gamma^{-1}\alpha\wedge\alpha\gamma+\gamma^{-1}
\alpha\wedge\extd\gamma+\gamma^{-1}\extd\gamma\wedge\gamma^{-1}
\alpha\gamma+\gamma^{-1}\extd\gamma\wedge\gamma^{-1}\extd\gamma}
using the Leibniz rule in the generalised form for $\Omega^2$ and the
assumption that $\wedge:\Omega^1\tens_M\Omega^1\to\Omega^2$. The 1st
and 8th,
3rd and 7th, and 4th and 9th terms cancel. The 6th term vanishes by
$\extd^2=0$. The point is that we only used the minimal structure for
$\Omega^1,\Omega^2$ as listed above.
The particular conventions here are read off from \cite[Fig.
3]{Ma:war95}.
Interestingly, by working with $A^{*\rm op}$ instead of $A^*$, all
the
computations can be done diagrammatically without any transpositions
or
braiding, i.e. the gauge theory at this level makes sense in any
monoidal
category with $\tens$ and $\oplus$ operations; see \cite{Ma:war95}.

We are now in a position to define the `fundamental group of a
differential
algebra' $(M,\Omega^1,\Omega^2)$. We fix
$M\equiv(M,\Omega^1,\Omega^2)$ as the
data which plays the role of a  differentiable manifold in the
formalism. Then
for all algebras $A^*$, we define
\eqn{flat}{ {\rm Flat}(M,A^*)=\{ \alpha\in\Omega^1\tens A^*| \
F(\alpha)=0\},}
the space of flat connections. The map $A^*\to {\rm Flat}(M,A^*)$
defines a
covariant functor from the category of algebras to the category of
sets.
Functors to sets are generally representable, which means there is
essentially
(up to some technical completions) an algebra $\pi_1(M)$ such that
\eqn{qhomot}{ \Hom_{\rm alg}(\pi_1(M),A^*)\isom{\rm
Flat(M,A^*)},\quad \forall
A^*}
This definition is based on the classical situation
\eqn{homot}{ \Hom_{\rm grp}(\pi_1(M),G)\isom{\rm Flat}(M,G),\quad
\forall G,}
where $G$ is a Lie group, $M$ is a manifold and we consider the
holonomy on
$P=M\times G$ defined by any flat connection.

This is the basic idea, and it has many variants. It is not known at
present
which variant works well or how to compute $\pi_1(M)$ in practice. As
explained
above, for a proper geometrical picture we need more structure on
$A^*$. For
example, a Hopf algebra structure and quantum tangent space for it.
Then (\ref{qhomot}) would be defined via $\Hom$ in the same category
and
$\pi_1(M)$ would likewise be a Hopf algebra equipped with a choice of
quantum
tangent space.

While clearly sketchy, we see from these remarks that one can in
principle use
quantum geometry to extract genuine `quantum topology' from a
noncommutative
algebra $M$. Although nontrivial computations are currently scarce,
the
$q$-monopole on $S_q^2$ constructed in \cite{BrzMa:gau} likewise
demonstrates
nontrivial topology, although from a slightly different point of view
(classically, the monopole charge is related to the fundamental group
of the
$S^1$ equator of $S^2$).

\subsection{Classical Integrable Systems}

In this section we specialise the above gauge theory further, to the
case where
$M$ is classical. In fact, we take $M=\C[\vecx]$, the coordinate ring
of
$\R^n$, but we allow non-standard $\Omega^1,\Omega^2$. Our goal is to
demonstrate the above theory. A second motivation comes from the
interesting
work of Muller-Hoissen\cite{DimMul:non}\cite{DMS:non}, where it is
already
shown that even this almost classical setup can have useful
applications. In
effect, we combine those ideas with the more systematic gauge theory
developed
in \cite{BrzMa:gau}. For our gauge group we take the classical group
with one
point $G=\{e\}$, so $A^*=\C$. This will make the point that
non-linearity
emerges not from non-Abelianness of the gauge group (though one can
have this
too), but from the breakdown of (\ref{comcalc}).

For our first example, we take $M=\C[x]$ 1-dimensional and $\Omega^1$
the 2-jet
calculus defined by $L=\{p,p^2\}$ as generated by $c(p)=p^3/3!$ in
Proposition~2.6. The relations of $\Omega^1$ are in (\ref{2-jet}).
Applying
$\extd$ and the Leibniz rule, we have in $\Omega^2$:
\eqn{2-jet2}{ \extd\omega=(\extd x)^2,\quad \omega^2=0,\quad
\omega\wedge\extd
x=-\extd x\wedge\omega,\quad x\extd\omega=(\extd\omega)x,\quad x\extd
x\wedge\omega=\extd x\wedge\omega x.}

\begin{propos} Working with the 2-jet differential calculus
(\ref{2-jet}),(\ref{2-jet2}), the curvature of a general gauge field
$\alpha=(\extd x)a(x) + \omega b(x)$ is
\[ F(\alpha)=\extd x\wedge\omega (b'-a''+2 a'\, a) - (\extd x)^2
(a'-b-a^2).\]
The gauge transformation by $\gamma(x)$ is
\[ a\to a+{\gamma'\over\gamma},\quad b\to
b-2a{\gamma'\over\gamma}+{\gamma''\over \gamma}
-2{(\gamma')^2\over\gamma^2},\]
under which $F$ is invariant. This is also the gauge symmetry of the
zero
curvature equation $a'=a^2+b$. The covariant derivative on a scalar
field
$\psi$ is
\[ \nabla\psi=\extd x(\psi'+a\psi)+\omega(\psi''+b\psi)\]
and is covariant under $\psi\to{\psi\over\gamma}$.
\end{propos}
\proof We write out $F=\extd \alpha +\alpha\wedge\alpha$ using the
Leibniz rule
and the relations (\ref{2-jet}),(\ref{2-jet2}) to collect terms. To
compute the
effect of gauge transformations, we use
(\ref{2-jet}) to find
\[ \gamma^{-1}(\extd x)a\gamma=(\extd x)a - 2\omega a
\gamma^{-1}\gamma'\]
\[ \gamma^{-1}\extd\gamma=\gamma^{-1}(\extd
x)\gamma'+\gamma^{-1}\omega\gamma''=(\extd
x)\gamma^{-1}\gamma'+\omega(\gamma^{-1}\gamma''-2\gamma^{-2}(\gamma')
^2).\]
It is a nontrivial check to verify that
$F(\alpha^\gamma)=\gamma^{-1}F\gamma =
F$ (the second equality by (\ref{2-jet2})), as it must by our general
theory.
The computation and covariance of $\nabla$ is equally simple. To
check
$\nabla^2\psi=F\psi$ one needs $\nabla$ on a general 1-form field
$\sigma=(\extd x)s+\omega t$, say. By similar computations to those
above, this
is
\eqn{covd1}{\nabla\sigma=(\extd x)^2(-s'+t+a^2)+\extd
x\wedge\omega(-s''+2a'a-bs+at+t').}
There are Bianchi identities for $F$ which one may check as well. All
of these
nontrivial facts are ensured by the gauge theory formalism in
Section~3.1.
\endproof

For our second example, we take $M=\C[x,y]$  and $\Omega^1$ the
$2$-dimensional
finite-difference calculus (\ref{2-dif}) on $\R^2$ defined by
$c(p,q)=\lambda^{-2}e^{\lambda p}+\mu^{-2}e^{\mu q}$. Applying
(\ref{2-dif})
and the Leibniz rule
to it, we obtain in $\Omega^2$ the identities $(\extd x)^2=0=(\extd
y)^2$ and
\eqn{2-dif2}{\extd x\wedge\extd y=-\extd
y\wedge\extd x,\quad x\extd x\wedge\extd y=(\extd x\wedge\extd
y)(x+\lambda),\quad y\extd x\wedge\extd y=(\extd x\wedge\extd
y)(y+\mu)}
Here, $\Omega^2$ is 1-dimensional and has a classical form.

\begin{propos} Working with the finite-difference calculus
(\ref{2-dif}),(\ref{2-dif2}), the curvature of a general gauge field
$\alpha=(\extd x)a(x,y)+(\extd y)b(x,y)$  is
\[ F(\alpha)=\extd x\wedge\extd
y\left(\del_{p_{1,0}}b-\del_{p_{0,1}}a+a(x,y+\mu)b -
b(x+\lambda,y)a\right).\]
The gauge transformation by $\gamma(x,y)$ is
\[a\to{a\gamma\over\gamma(x+\lambda,y)}+{1\over\lambda}
\left(1-{\gamma\over\gamma(x+\lambda,y)}\right),\quad
b\to{b\gamma\over\gamma(x,y+\mu)}+{1\over\mu}
\left(1-{\gamma\over\gamma(x,y+\mu)}\right),\]
under which $F\equiv(\extd x\wedge\extd y)F(x,y)$ transforms as
\[ F\to F{\gamma\over\gamma(x+\lambda,y+\mu)}.\]
This is also the gauge symmetry of the zero curvature equation. The
covariant
derivative on a scalar field $\psi(x,y)$ is
\[ \nabla\psi=\extd x (\del_{p_{1,0}}\psi +a\psi)+\extd
y(\del_{p_{0,1}}\psi+b\psi)\]
and is covariant under $\psi\to {\psi\over\gamma}$.
\end{propos}
\proof We use the relations $f(x,y)\extd x=(\extd x) f(x+\lambda,y)$
etc. to
find the coefficient of $\extd x\wedge\extd y$ in
$\extd\alpha+\alpha\wedge\alpha$. This also gives
$\gamma^{-1}\alpha\gamma$ and
$\gamma^{-1}F\gamma$. Finally, to compute the remaining part of the
gauge
transformation, we have
\[ \gamma^{-1}\extd\gamma=\gamma^{-1}(\extd
x){\gamma(x+\lambda,y)-\gamma(x,y)\over\lambda}+\gamma^{-1}(\extd
y){\gamma(x,y+\mu)-\gamma(x,y)\over\mu},\]
using (\ref{2-dif}). Moving $\gamma^{-1}$ to the right then gives the
form
shown. Again, it is a nontrivial check on the calculus that $F$ is
indeed
covariant under this gauge transformation. The covariant derivative
on a 1-form
field $\sigma=(\extd x)s(x,y)+(\extd y)t(x,y)$ is likewise computed,
as
\[ \nabla\sigma=\extd x\wedge\extd
y\left(\del_{p_{1,0}}t-\del_{p_{0,1}}s+a(x,y+\mu)t
-b(x+\lambda,y)s\right),\]
from which one may verify that $\nabla^2\psi=F\psi$. These and the
other
properties are all ensured by the gauge theory formalism in
Section~3.1.
\endproof

Although these examples are elementary (and Proposition~3.2 is
probably not new
to anyone who has played with lattice gauge theory), we see that
non-linear
differential equations and finite-difference equations on $\R^n$ can
be treated
as actual zero-curvature equations for a general formalism of gauge
theory;
just with a choice of non-standard calculi. As such, assuming trivial
`quantum
fundamental group' $\pi_1(M)$ (as defined in Section~3.1), they are
completely
integrable in the
sense that every solution would be of the form
$\alpha=\gamma^{-1}\extd\gamma$
as the gauge transform of the zero solution. It is not clear what
class of
integrable systems can be covered in this way, but by using the
higher
jet-calculi
they can include arbitrary derivatives. By using a more general
$c(p)$ we have entirely new gauge theories as well. The point is that
one does
not need to invent and verify each of these generalised gauge
theories `by
hand'; they are instead constructed as examples of a single
formalism\cite{BrzMa:gau}\cite{BrzMa:coa}\cite{Ma:diag}.

Moreover, the formalism allows for the possibility of non-trivial
gauge
symmetry, which can be a group, quantum group, braided group or even
a general
algebra. For example, a complete formalism of lattice non-Abelian
gauge theory
and scalar, vector matter fields is an immediate corollary. And we
are not
limited to classical spaces $M$ for the base. They can be quantum,
super or
even anyonic (where the coordinate obeys $\theta^n=0$). An example of
braided
gauge theory on an anyonic base is worked out in detail in
\cite{Ma:diag}.
Finally, we are not limited to trivial bundles. Bundles can be
nontrivial both
in a familiar geometrical way and in a purely quantum way (even when
they are
geometrically trivial), the latter being controlled  by a certain
nonAbelian
2-cohomology; \cite{Ma:war95} is a recent review.

\subsection*{Acknowledgements} It is a pleasure to thank H.-D.
Doebner, P.
Kramer, W. Scherer, V. Dobrev and other members of the organising
committee for
an excellent conference Group XXI in Goslar.


 \end{document}